\documentclass[11pt]{amsart}
\usepackage{amsfonts}
\usepackage{epsfig,verbatim}

\theoremstyle{plain}
\newtheorem{theorem}{Theorem} [section]
\newtheorem{proposition} [theorem] {Proposition}
\newtheorem{corollary} [theorem] {Corollary}
\newtheorem{lemma} [theorem] {Lemma}

\theoremstyle{definition}

\newtheorem{definition} [theorem] {Definition}

\theoremstyle{remark}
\newtheorem{remark} [theorem] {Remark}

\def \r{\mbox{${\mathbb R}$}}

\def \h{\mbox{${\mathbb H}$}}
\def \e{\mbox{${\mathbb E}$}}
\def \s{\mbox{${\mathbb S}$}}

\DeclareMathOperator{\trace}{trace}
\DeclareMathOperator{\cst}{constant}

\begin{document}

\begin{abstract}
In this article we characterize all biharmonic curves of the
Cartan-Vranceanu $3$-dimensional spaces and we give their explicit parametrizations.
\end{abstract}

\title[Biharmonic curves]{The classification of biharmonic curves of\\  Cartan-Vranceanu $3$-dimensional spaces}

\author{R.~Caddeo}
\author{S.~Montaldo}
\author{C.~Oniciuc}
\author{P.~Piu}

\address{Universit\`a degli Studi di Cagliari\\
Dipartimento di Matematica e Informatica\\
Via Ospedale 72\\
09124 Cagliari}
\email{caddeo@unica.it,montaldo@unica.it,piu@unica.it}

\address{Faculty of Mathematics\\ ``Al.I. Cuza'' University of Iasi\\
Bd. Carol I no. 11 \\
700506 Iasi, ROMANIA}
\email{oniciucc@uaic.ro}

\thanks{Work partially supported by GNSAGA (ITALY); the third author was supported by a
CNR-NATO fellowship (ITALY), and by the Grant At, 73/2005, CNCSIS (ROMANIA)}

\subjclass{58E20}
\keywords{Harmonic and biharmonic maps, helices}

\maketitle

\section{Introduction}

Biharmonic curves $\gamma:I\subset\r\to (N,h)$ of a
Riemannian manifold are the solutions of the fourth
order differential equation

$$
\nabla^3_{\gamma'}\gamma' - R(\gamma',\nabla_{\gamma'} \gamma')\gamma' = 0.
$$

\noindent As we shall detail in the next  section,
they arise from a variational problem and are a
natural generalization of geodesics.

    In the last decade have appeared several
    papers on the construction and classification
    of biharmonic curves starting with \cite{RCSMPP1},
    where the authors described the case of curves of a surface.
    Biharmonic curves in a $3$-dimensional Riemannian manifold with constant sectional
    curvature $K\leq 0$ are geodesics (see \cite{Dim},
     for $K=0$, and \cite{RCSMCO2}, for $K<0$);
while, in \cite{RCSMCO1}, the authors proved that for $K>0$
the biharmonic curves  are helices, that is curves with
constant geodesic curvature and geodesic torsion.

    Among the $3$-dimensional manifolds of non-constant
    sectional curvature
a special role is played by the homogeneous Riemannian spaces
with a large isometry group. For these spaces, except for
those with constant negative curvature, there is a nice
local representation given by the following two-parameter
family of Riemannian metrics
(the Cartan-Vranceanu metric)
$$
 ds^{2}_{\ell,m} =\frac{dx^{2} + dy^{2}}{[1 + m(x^{2} + y^{2})]^{2}} +  \left(dz +
\frac{\ell}{2} \frac{ydx - xdy}{[1 + m(x^{2} + y^{2})]}\right)^{2},\quad \ell,m \in {\r}
$$
\noindent defined on $M=\r^3$ if $m\geq 0$, and on
$M=\{(x,y,z)\in\r^3\;:\; x^{2} + y^{2} < - \frac{1}{m}\}$
otherwise.

\noindent Biharmonic curves on $(M,ds^{2}_{\ell,m})$ have
been already studied for particular values of $\ell$ and $m$.
In particular, if $ m = 0$ and $\ell \neq 0$, $(M,ds^2 _{\ell,m})$
is the Heisenberg space $\h_3$ endowed with a left invariant
metric and the explicit solutions of the biharmonic curves
were obtained in \cite{COP};
if $\ell = 1$ and $m\neq 0$ a study of the biharmonic curves
was given also in \cite{Inoguchi}.

    In this paper we prove that the biharmonic curves of
    $(M,ds^{2}_{\ell,m})$
are helices and we find out their explicit parametric
equations, for all values of $\ell$ and $m$.

\section{Preliminary}

\subsection{Biharmonic curves}
Harmonic maps $\phi : (M,g)\to (N,h)$ between
Riemannian manifolds  are the critical points of the energy functional
$E(\phi)=\frac{1}{2}\int_{M}\, |d\phi|^2\,v_g$, and the  corresponding
Euler-Lagrange equation is given by the vanishing of the tension
field $\tau(\phi)=\trace\nabla d\phi$.
Biharmonic maps (as suggested by J.~Eells and J.H.~Sampson
in \cite{ES}) are the critical points of the {\it bienergy} functional
$E_2(\phi)=\frac{1}{2}\int_{M}\, |\tau(\phi)|^2\,v_g$.
In \cite{Jia} G.Y.~Jiang derived the first variation formula
of the bienergy showing that the Euler-Lagrange equation
for $E_2$ is
\begin{eqnarray*}
\tau_2(\phi) &=& - J(\tau(\phi))=
-\Delta\tau(\phi) -\trace R^{N}(d\phi,\tau(\phi))d\phi\notag \\
&=&0,
\end{eqnarray*}
where $J$ is the Jacobi operator of $\phi$ and
$R^{N}(X,Y)=\nabla_X\nabla_Y-\nabla_Y\nabla_X-\nabla_{[X,Y]}$.
The equation $\tau_{2}(\phi)=0$ is called the
{\it biharmonic equation}.

\noindent Since $J$ is linear, any harmonic map is biharmonic.
Therefore, the main interest is to find and classify {\it proper}
biharmonic maps, that is non-harmonic biharmonic maps.

    In this paper we restrict our attention to  curves
$\gamma:I\to (N,h)$  parametrized by arc length, from an open
interval $I\subset\r$ to a Riemannian
manifold. In this case, putting $T=\gamma'$, the tension field
becomes
$ \tau({\gamma})=\nabla_{T}T$ and the
biharmonic equation reduces to
\begin{equation}\label{bieq-curve}
\nabla^3_{T}T - R(T,\nabla_T T)T = 0.
\end{equation}

\noindent To describe geometrically Equation \eqref{bieq-curve}
let recall the definition of the Frenet frame.
\begin{definition}[See, for example, \cite{Laugwitz}] \label{def2.1}
The Frenet frame $\{F_{i}\}_{i=1,\dots,n}$ associated to a
curve $\gamma : I\subset {\r}\to (N^{n},h)$  parametrized
by arc length is the orthonormalisation of the $(n+1)$-uple
$\{ \nabla_{\frac{\partial}{\partial t}}^{(k)} d\gamma
(\frac{\partial}{\partial t}) \}_{k=0,\dots,n}$, described
by:
\begin{align*} F_{1}&=d\gamma
(\frac{\partial}{\partial t}) ,  \\
\nabla_{\frac{\partial}{\partial
t}}^{\gamma} F_{1} &= k_{1} F_{2} , \\
\nabla_{\frac{\partial}{\partial
t}}^{\gamma} F_{i} &= - k_{i-1} F_{i-1}
+ k_{i}F_{i+1} , \quad \forall i =
2,\dots,n-1 , \\
\nabla_{\frac{\partial}{\partial
t}}^{\gamma} F_{n} &= - k_{n-1} F_{n-1}
, \end{align*}
where the functions $\{k_{1},k_{2},\ldots,k_{n-1}\}$
are called the curvatures of $\gamma$ and $\nabla^{\gamma}$
is the connection on the pull-back bundle $\gamma^{-1}(TN)$. Note that
$F_{1}=T={\gamma}'$ is the unit tangent vector field along the curve.
 \end{definition}

\noindent Using the Frenet frame, the biharmonic equation
\eqref{bieq-curve} reduces to a differential system of the
curvatures of $\gamma$ as shown in the following

\begin{proposition}\label{prop1}
Let $\gamma : I \subset {\r} \to (N^{n},h)$ ($n\geq 2$) be
a curve parametrized by arc length from an open interval
of $\r$ into a Riemannian manifold $(N,g)$.
Then $\gamma$ is biharmonic if and only if:
$$
\left\{
\begin{array}{l}
{k_{1}} {k'_{1}} = 0\\
{k''_1}-k_{1}^3 - k_{1} k_2^2  +  k_{1} R(F_1,F_2,F_1,F_2) = 0 \\
2 {k'_{1}} k_2 + k_1 {k'_2} + k_1 R(F_1,F_2,F_1,F_3)  = 0\\
k_{1} k_2 k_3 + k_1 R(F_1,F_2,F_1,F_4) = 0 \\
k_1 R(F_1,F_2,F_1,F_j) = 0 \hspace{1,5 cm} j = 5,\ldots, n
\end{array}
\right.
$$
\end{proposition}
\begin{proof}
With respect to its Frenet frame, the biharmonic
equation of $\gamma$ is:
\begin{eqnarray*}
\nabla^3_{F_1}F_1 - R(F_1,\nabla_{F_1} F_1)F_1 &=&-3 k_{1} {k'_{1}} F_1 +
({k''_1} - k_{1}^3 - k_{1} k_2^2 )F_2 \nonumber\\
&&+(2 {k'_{1}} k_2 + k_1 {k'_2})F_3 + k_{1} k_2 k_3 F_4 - k_1 R(F_1,F_2)F_1\\
&=& 0. \nonumber
\end{eqnarray*}

\end{proof}

\noindent If we look for proper biharmonic solutions,
that is for biharmonic
curves with $k_{1}\neq0$, we have
\begin{equation}\label{euler-lagrange}
\left\{
\begin{array}{l}
k_{1} = \cst\neq 0\\
k_{1}^2 + k_2^2 =  R(F_1,F_2,F_1,F_2) \\
{k'_2} = -R(F_1,F_2,F_1,F_3) \\
k_2 k_3 = - R(F_1,F_2,F_1,F_4) \\
R(F_1,F_2,F_1,F_j) = 0 \hspace{1,5 cm} j = 5,\ldots, n
\end{array}
\right.
\end{equation}

\subsection{Riemannian structure of Cartan-Vranceanu $3$-manifolds}
Let $m$ be a real parameter. We shall denote by $M$
the whole $\r^3$ if $m\geq 0$, and by
$M=\{(x,y,z)\in\r^3\;:\; x^{2} + y^{2} < - \frac{1}{m}\}$  otherwise.
Consider on $M$ the following two-parameter
family of Riemannian metrics

\begin{equation}\label{1.1}
 ds^{2}_{\ell,m} =\frac{dx^{2} + dy^{2}}{[1 + m(x^{2} + y^{2})]^{2}} +  \left(dz +
\frac{\ell}{2} \frac{ydx - xdy}{[1 + m(x^{2} + y^{2})]}\right)^{2},
\end{equation}
where $\ell,m \in {\r}$.

\noindent These metrics have been known for a long time.
They can be found in the
classification of $3$-dimensional homogeneous metrics
given by L. Bianchi in
1897 (see \cite{Bi}); later, they appeared in form (\ref{1.1})
in  \'{E}. Cartan, (\cite{Ca}
p. 304) and in G. Vranceanu (see \cite{Vr}, p. 354).
Their geometric interest lies in the
following fact: {\it the family of metrics (\ref{1.1})
includes all $3$-dimensional
homogeneous metrics whose group of isometries has dimension $4$ or
$6$, except for those of constant negative sectional curvature}.

\noindent The Cartan-Vranceanu metric \eqref{1.1} can be written as:
$$
 ds^{2}_{\ell,m} =\sum_{i=1}^3 \omega^i \otimes \omega^i
$$
where, putting $F=1 + m(x^{2} + y^{2})$,
\begin{equation}\label{forme}
\omega^1 = \frac{dx}{F}, \quad \omega^2 = \frac{dy}{F},
\quad \omega^3 = dz + \frac{\ell}{2} \frac{ydx - xdy}{F},
\end{equation}
and the orthonormal basis of dual vector fields to the
$1$-forms \eqref{forme} is
\begin{equation}\label{baseorto}
E_1 = F \frac{\partial}{\partial x} - \frac{\ell\,y}{2}
\frac{\partial}{\partial z},\quad E_2 = F
\frac{\partial}{\partial y} +\frac{\ell\,x}{2}
\frac{\partial}{\partial z}, \quad E_3 =
\frac{\partial}{\partial z}\,.
\end{equation}

For completeness we give the expressions, with respect to
the orthonormal basis \eqref{baseorto}, of the Levi-Civita
connection and of the nonzero components of the curvature
and Ricci tensors:
\newpage
\begin{equation}\label{eq:conn}
\begin{array}{l|l|l}
\nabla_{E_1}E_1 = 2 m y E_2 & \nabla_{E_1}E_2=- 2 m y E_1 + \frac{\ell}{2}E_3 &\nabla_{E_1}E_3=\nabla_{E_3}E_1= - \frac{\ell}{2}E_2\\
\nabla_{E_2}E_2 = 2 m x E_1 & &
\\
\nabla_{E_3}E_3 = 0 & \nabla_{E_2}E_1= - 2m x E_2 - \frac{\ell}{2}E_3   & \nabla_{E_2}E_3 =  \nabla_{E_3}E_2=\frac{\ell}{2}E_1
\end{array}
\end{equation}
\begin{equation}\label{curvature}
R_{1212} = 4 m - \frac{3}{4}\ell^2, \quad R_{1313} =
\frac{\ell^2}{4}, \quad R_{2323} =  \frac{\ell^2}{4}
\end{equation}
\begin{equation}\label{ricci}
\rho_{11} = \rho_{22} = 4m - \frac{\ell^2}{2}, \quad
\rho_{33} = \frac{\ell^2}{2}.
\end{equation}
\begin{remark}\label{re-cv} \quad\\
\begin{itemize}
\item If $\ell=0$, then $M$ is the product of a surface
$S$  with constant Gaussian curvature $4m$ and the real line $\r$.
\item If $\rho_{11} - \rho_{33} = 4m - \ell^2 = 0 $, then
$M$ has non negative constant sectional curvature.
\item From the Kowalski's classification \cite{OK1} we
know that the principal Ricci curvatures of $SU(2)$
satisfy
$\rho_{33} > 0$, $\rho_{11} + \rho_{33} > 0$ and  $\rho_{11} \neq \rho_{33}$.
Thus, form \eqref{ricci}, if $\ell \neq 0$ and $m>0$, $M$ is
locally  $SU(2)$.
\item Similarly, if  $\ell \neq 0$ and $m < 0$, $M$ is
locally $\widetilde{SL}(2,\r)$, while if $ m = 0$ and $\ell \neq 0$
we get the left invariant metric on the Heisenberg space $\h_3$.
\end{itemize}
\end{remark}

\begin{figure}[htb]
\centerline{
\begin{picture}(100,160)(-40,-100)
\put(0,0){\qbezier(-60,25)(0,-100)(60,25)}
\put(0,-70){\vector(0,1){100}}
\put(-70,-37.5){\vector(1,0){140}}
\put(3,22){\makebox{$m$}}
\put(65,-35){\makebox{$\ell$}}
\put(72,-38){\makebox{$\h_{3}$}}
\put(-12,32){\makebox{$\s^{2}\times\r$}}
\put(-14,-80){\makebox{$\h^{2}\times\r$}}
\put(62,25){\makebox{$\s^{3}$}}
\put(-75,26){\makebox{$4m = \ell^2$}}
\put(2,-35){\makebox{$\r^{3}$}}
\put(-70,-60){\makebox{${\widetilde{SL}(2,{\r})}$}}
\put(30,-60){\makebox{${\widetilde{SL}(2,{\r})}$}}
\put(-70,-20){\makebox{$SU(2)$}}
\put(40,-20){\makebox{$SU(2)$}}
\put(10,10){\makebox{$SU(2)$}}
\put(-40,10){\makebox{$SU(2)$}}
\put(0,-37.5){\circle*{4}}
\end{picture}
}
\caption{\label{diagram} The geometric description of
the metric $ds^2_{\ell,m}$.}
\end{figure}
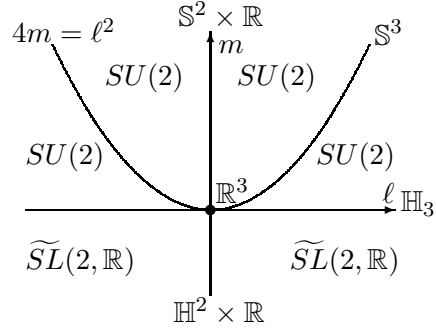

\noindent Remark~\ref{re-cv} gives rise to a  nice geometric
description of the metric $ds^2_{\ell,m}$ as shown in Figure~\ref{diagram}.

\section{Biharmonicity conditions for curves in $(M,ds^2 _{\ell,m})$}
Let $\gamma:I\to (M,ds^2 _{\ell,m})$ be a differentiable curve
parametrized by arc length and let $\{F_1=T = T_iE_i,\,F_2=N = N_iE_i,\,
F_3=B =  B_iE_i\}$ be the Frenet frame
field tangent to $M$ along $\gamma$ decomposed with
respect to the orthonormal basis \eqref{baseorto}.

\noindent By making use of \eqref{euler-lagrange}
 and of the expression  of the
curvature tensor field \eqref{curvature}, we obtain the
 following system for the proper biharmonic curves
\begin{equation}\label{3.3Vran}
\left\{
\begin{array}{l}
k=\cst\neq 0 \\
k^2+\tau^2=\frac{\ell^2}{4}-(\ell^2 - 4m) B_3^2 \\
\tau'= (\ell^2 - 4m)N_3B_3 ,
\end{array}
\right.
\end{equation}
where $k=k_1$ and $\tau=-k_2$.

    By analogy with  curves in $\r^3$, also
    following \cite{GL},
we keep
the name {\it helix}
for a curve in a Riemannian manifold having constant
both geodesic curvature
and geodesic torsion.

\begin{remark}\label{re:known} \quad
\begin{itemize}
\item[({\it i})] If $\ell = m = 0$, $(M,ds^2 _{\ell,m})$
is the Euclidean space and $\gamma$ is biharmonic if and
only if it is a line (see \cite{Dim});
\item[({\it ii})] if $\ell^2 = 4m$ and $\ell \neq 0$,
then $(M,ds^2 _{\ell,m})$ is locally the $3$-dimensional
sphere and the proper biharmonic curves were classified
in  \cite{RCSMCO1}, where it was proved that they are helices;
\item[({\it iii})] if $\ell = 0$ and $m < 0$, $(M,ds^2 _{\ell,m})$
is isometric to $\h^2 \times\r$ with the product metric and it can be show that  all
biharmonic curves are geodesics.
\item[({\it iv})] if $ m = 0$ and $\ell \neq 0$, $(M,ds^2 _{\ell,m})$
is the Heisenberg space $\h_3$ endowed with a left
invariant metric and the biharmonic curves were studied
in \cite{COP}.
\item[({\it v})] if $\ell = 1$ a study of the explicit
solutions of \eqref{3.3Vran} was given in \cite{Inoguchi}.
\end{itemize}
\end{remark}

    From now on we shall assume that $\ell^2\neq 4m$ and $m\neq 0$.
    This is, essentially, the only case left to study according to Remark~\ref{re:known}.

\noindent As in previous cases we have

\begin{theorem}\label{gener-prop}
If $\gamma:I\to (M,ds^2 _{\ell,m})$ is a proper biharmonic
curve parametrized by arc length, then it is a helix.
\end{theorem}

\begin{proof}
Let  $\gamma:I\to (M,ds^2 _{\ell,m})$ be a non geodesic
curve parametrized by arc length. Then from the Frenet's equation we have
$$
\langle\nabla_T B,E_3\rangle = \tau N_3\,,
$$
while, using the definition of covariant derivative, we get
$$
\langle\nabla_T B,E_3\rangle = B'_3 + \frac{\ell}{2}(T_1 B_2 - T_2 B_1) = {B'}_3 - \frac{\ell}{2} N_3.
$$
Comparing the two equations we have
\begin{equation}\label{eq-taun3}
\tau N_3 = {B'}_3 - \frac{\ell}{2} N_3\,.
\end{equation}

Assume now that $\gamma$ is a proper biharmonic curve. We
first show that  $B_3 \neq 0$. Indeed, if $B_3 = 0$, then
\eqref{eq-taun3} implies that
$$
 N_3 (\tau + \frac{\ell}{2} ) = 0.
$$
The latter equation gives two possibilities:
\begin{itemize}
\item if $N_3=0$ (and $B_3=0$) then  $T = \pm E_3$ and the
curve $\gamma$ is a geodesic;
\item if $(\tau + \frac{\ell}{2} ) = 0$ then, using the
second equation of \eqref{3.3Vran}, we must have that
$\gamma$ is again a geodesic.
\end{itemize}

\noindent Therefore, $B_3 \neq 0$. Deriving the second
equation of \eqref{3.3Vran} yields
$$
\tau \tau' = - (\ell^2 - 4 m) B_3 B'_3
$$
and, taking into account the third equation of \eqref{3.3Vran}, gives
$$
(\ell^2 - 4 m) B_3 (\tau N_3 + B'_3) = 0.
$$
Thus
$
\tau N_3 = - B'_3
$
that summed with
$
\tau N_3 =  B'_3 - \frac{\ell}{2} N_3
$
leads to
$$
N_3 (4 \tau + \ell) = 0,
$$
and consequently $\tau$ is constant.

\end{proof}

\noindent From the proof of Theorem~\ref{gener-prop}
and \eqref{3.3Vran}  we
have, in conclusion,

\begin{corollary}
\label{eq:biharmoniccurves} Let $\gamma:I\to (M,ds^2_{\ell,m})$
be a curve parametrized by arc length.
Then $\gamma$ is a proper biharmonic curve if
and only if
\begin{equation}\label{eq:biharmonichelix}
\left\{
\begin{array}{l}
k = \cst\neq 0 \\
\tau = \cst \\
N_3 = 0 \\
k^2+\tau^2=\frac{\ell^2}{4}-(\ell^2 - 4 m)B_3^2 .
\end{array}
\right.
\end{equation}
\end{corollary}

\section{Explicit formulas for proper biharmonic
curves in $(M,ds^2_{\ell,m})$}

In this section we use Corollary~\ref{eq:biharmoniccurves}
to derive the explicit parametric
equations of proper biharmonic curves in  $(M,ds^2_{\ell,m})$.
We first prove the following
\begin{lemma}\label{lemma1}
Let $\gamma:I\to (M,ds^2_{\ell,m})$ be a non-geodesic curve
parametrized by arc length.
If $N_3=0$, then
\begin{equation}\label{tangelica}
T(t) =\sin\alpha_0\cos\beta(t) E_1+\sin\alpha_0\sin\beta(t)
E_2+\cos\alpha_0 E_3,
\end{equation}
where $\alpha_0\in (0,\pi)$.
\end{lemma}

\begin{proof}
If $\gamma'=T=T_1E_1+T_2E_2+T_3E_3$,
from
\begin{eqnarray*}
\nabla_TT &=& (T_1'+ \ell\, T_2\,T_3 + 2\, m\, x\, T^2_2 - 2\, m\, y\, T_1\, T_2)\, E_1 \\
&& +(T_2' - \ell\, T_1\,T_3 + 2\, m\, y\, T^2_1
-2\, m\, x\, T_1\, T_2 ) E_2 + T_3'\, E_3 \\
&=& kN
\end{eqnarray*}
it follows that $N_3=0$ if and only if $T_3' = 0$, i.e. if and only if $T_3
=\cst$.
Since $T_3\in (-1,1)$ and the norm of $T$ is one, there
exists a constant $\alpha_0\in(0,\pi)$ and a unique
(up to a additive constant $2k \pi$) smooth function $\beta$ such that
$$
T(t) =\sin\alpha_0\cos\beta(t) E_1+\sin\alpha_0\sin\beta(t)
E_2+\cos\alpha_0 E_3
$$
\end{proof}
Starting from the expression \eqref{tangelica} of $T$ we
are ready to state the main result.
\begin{theorem}
\label{eq:biharmonicintegralcurves}
Let $(M,ds^2_{\ell,m})$ be the Cartan-Vranceanu space
with $m\neq 0$ and
$\ell^2-4m\neq 0$.
Assume that
$\delta=\ell^2+(16 m - 5 \ell^2)\sin^2\alpha_0\geq 0$,
$\alpha_0\in(0,\pi)$, and denote  by $2\omega_{1,2}=-\ell \cos\alpha_0\pm \sqrt{\delta}$.
Then, the  parametric equations of all proper biharmonic curves of
$(M,ds^2_{\ell,m})$ are of the following three types.

\noindent {\bf Type I}
\begin{equation}
\label{eq:biharmonichelices1}
\left\{
\begin{array}{l}
x(t) = \displaystyle{b \sin \alpha_0 \sin \beta(t) + c}, \quad b,c\in\r,\,b>0\\
y(t) =\displaystyle{ - b \sin \alpha_0 \cos \beta(t)+d}, \quad d\in\r\\
z(t) = \displaystyle{\frac{\ell}{4m} \beta(t)+ \frac{1}{4m}\left[ (4m - \ell^2)\cos \alpha_0  - \ell\, \omega_{1,2}\right]t },
\end{array}
\right.
\end{equation}
where
$\beta$ is a non-constant solution of the following ODE:
\begin{equation}\label{eq-diff-beta}
\beta'+2\, m\, d\, \sin\alpha_0 \cos\beta - 2\, m\, c\, \sin\alpha_0 \sin\beta=\ell\, \cos\alpha_0 +
2\, m\, b\, \sin^2\alpha_0 + \omega_{1,2},
\end{equation}
and the constants satisfy
$$
c^2+d^2=\frac{b}{m}\left\{(\ell \cos\alpha_0+\omega_{1,2}-\frac{1}{b})+m\,b \sin^2\alpha_0\right\}.
$$
\noindent {\bf Type II}
If $\beta=\beta_0=\cst$ and $\cos\beta_0 \sin\beta_0\neq 0$,
the parametric equations are
$$
\label{eq:biharmonichelicesII}
\left\{
\begin{array}{l}
x(t) = \displaystyle{x(t)}\\
y(t) =\displaystyle{ x(t) \tan\beta_0 +a} \\
z(t) = \displaystyle{\frac{1}{4m}\left[ (4m - \ell^2)\cos \alpha_0  - \ell\, \omega_{1,2}\right]t +b},\quad b\in\r
\end{array}
\right.
$$
where $a=\frac{\omega_{1,2}+\ell \cos\alpha_0}{2m \sin\alpha_0 \cos\beta_0}$
and $x(t)$ is a solution
of the following ODE:
\begin{equation}\label{eq-diff-betaII}
x'=\big(1+m[x^2+(x \tan\beta_0+a)^2]\big)\sin\alpha_0 \cos\beta_0\,.
\end{equation}

\noindent {\bf Type III} If $\cos\beta_0 \sin\beta_0 = 0$,
up to interchange of $x$ with $y$, $\cos\beta_0=0$ and the
parametric equations are
$$
\label{eq:biharmonichelicesIII}
\left\{
\begin{array}{l}
x(t) = x_0= \mp\displaystyle{ \frac{\omega_{1,2}+\ell \cos\alpha_0}{2 m \sin\alpha_0}} \\
y(t)= y(t)\\
z(t) = \displaystyle{\frac{1}{4m}\left[ (4m - \ell^2)\cos \alpha_0  - \ell\, \omega_{1,2}\right]t +b},\quad b\in\r
\end{array}
\right.
$$
where $y(t)$ is a solution
of the following ODE:
\begin{equation}\label{eq-diff-betaIII}
y'=\pm\big(1+m[x_0^2+y^2]\big)\sin\alpha_0\,.
\end{equation}
\end{theorem}
\begin{proof}
Let $\gamma(t)=(x(t),y(t),z(t))$ be a biharmonic curve
parametrized by arc length. We shall make use of the Frenet
formulas, and we shall take into account
Corollary~\ref{eq:biharmoniccurves} and Lemma~\ref{lemma1}.
The covariant derivative of the vector field $T$ given by
\eqref{tangelica} is
\begin{eqnarray*}
\nabla_TT
&=&\big[- \beta' \sin \alpha_0 \sin\beta - 2 m y \sin^2 \alpha_0 \cos\beta\sin\beta \\
& &+
2 m x\sin^2 \alpha_0 \sin^2 \beta + \ell \cos \alpha_0 \sin \alpha_0 \sin\beta\big] E_1 \\
& &+ \big[\beta' \sin \alpha_0 \cos\beta + 2 m y \sin^2 \alpha_0 \cos^2 \beta -\\
& & - 2 m x\sin^2 \alpha_0 \cos\beta \sin \beta - \ell \cos \alpha_0 \sin \alpha_0 \cos\beta\big] E_2 \\
&=& kN,\\
\end{eqnarray*}
where
$$
k=\vert \beta'   + 2 m y \sin \alpha_0 \cos \beta -
2 m x\sin \alpha_0 \sin \beta  - \ell \cos \alpha_0 \vert \sin \alpha_0.
$$
We assume that
\begin{equation}\label{eq-omega-ex}
\omega=\beta'  + 2 m y \sin \alpha_0 \cos \beta -
2 m x\sin \alpha_0 \sin \beta  - \ell \cos \alpha_0  > 0.
\end{equation}
Then we have
\begin{equation}
\label{eq:curburag}
k = \omega \sin \alpha_0
\end{equation}
and
$$
N=- \sin \beta E_1 + \cos \beta E_2.
$$
Next,
\begin{equation}\label{eq:binorm}
B = T\times N = - \cos\alpha_0  \cos \beta E_1 - \cos \alpha_0 \sin \beta E_2 +
\sin \alpha_0 E_3
\end{equation}
and

\begin{eqnarray*}
\nabla_T B &=& [\beta' \cos \alpha_0 \sin \beta  + 2 m y \sin \alpha_0 \cos \alpha_0 \sin \beta \cos \beta -
2 m x\cos \alpha_0 \sin \alpha_0 \sin^2 \beta  \\
&&-\frac{\ell}{2} \cos^2 \alpha_0 \sin \beta + \frac{\ell}{2} \sin^2 \alpha_0 \sin\beta]E_1 \\
&&+
[-\beta' \cos \alpha_0 \cos \beta  - 2 m y \sin \alpha_0 \cos \alpha_0 \cos^2 \beta  +
2 m x \cos \alpha_0 \sin \alpha_0 \sin\beta \cos\beta  \\
&&+ \frac{\ell}{2} \cos^2 \alpha_0 \cos\beta - \frac{\ell}{2} \sin^2 \alpha_0 \cos\beta]E_2\,.
\end{eqnarray*}
It follows that the geodesic torsion $\tau$ of $\gamma$ is given by
\begin{equation}
\label{eq:torsiuneag}
\tau = - \omega \cos \alpha_0 - \frac{\ell}{2}.
\end{equation}
In order to find the explicit equations for
$\gamma(t) = (x(t),y(t),z(t))$,
we must integrate the system
${d\gamma}/{dt}=T$, that in our case is
\begin{equation}\label{eq.bihdif}
\begin{cases}
\displaystyle{\frac{{x}'}{1 + m(x^2 + y^2)}=\sin\alpha_0\cos\beta}
 \\
\displaystyle{\frac{{y}'}{1 + m(x^2 + y^2)}=\sin\alpha_0\sin \beta}
 \\
\displaystyle{{z}'=\cos\alpha_0
+\frac{\ell}{2}\sin\alpha_0\big(x \sin\beta-y \cos \beta\big )}.
\end{cases}
\end{equation}
We now assume that $\beta'\neq 0$, that is we considerer solutions
of
Type I.
Deriveting \eqref{eq:curburag} and taking into account \eqref{eq.bihdif}
we get
$$
\beta'' = \beta'  \frac{ 2mx {x}' +2my {y}'}{1 + m (x^2 + y^2)}.
$$
By integration of last equation we find
\begin{equation}\label{eq-b}
 1 + m (x^2 + y^2) = b \beta'\,,\quad b>0.
\end{equation}
Replacing \eqref{eq-b} in  \eqref{eq.bihdif} and integrating
we obtain the  solution
\begin{equation}
\label{eq:biharmonichelices2}
\left\{
\begin{array}{l}
x(t) = \displaystyle{b \sin \alpha_0 \sin \beta(t) + c} \\
y(t) =\displaystyle{ - b \sin \alpha_0 \cos \beta(t)+d} \\
z(t) = \displaystyle{(\cos \alpha_0 +\frac{\ell\,b}{2}\sin^2 \alpha_0)t
+\frac{\ell}{2} \int\left(c \sin\alpha_0\sin\beta-d\sin\alpha_0\cos\beta\right)\, dt}\,.
\end{array}
\right.
\end{equation}
To determine $\beta$ we replace in
\eqref{eq:biharmonichelix} the values of $k$, $\tau$ and $B_3$
given in  \eqref{eq:curburag}, \eqref{eq:torsiuneag} and
\eqref{eq:binorm} respectively. This gives
\begin{equation}\label{eq-omega}
\omega^2+\omega\,\ell\,\cos\alpha_0+(\ell^2-4m)\sin^2\alpha_0=0
\end{equation}
Assume that
$\delta=\ell^2+(16 m - 5 \ell^2)\sin^2\alpha_0\geq 0$; then
the solutions of
\eqref{eq-omega} are
$$
\omega_{1,2}=\frac{-\ell \cos\alpha_0\pm \sqrt{\delta}}{2},
$$
which are always different from zero.
Since we have assumed that $\omega$ is a positive constant,
we have to choose
the positive root of \eqref{eq-omega}. We point out that
if $\omega$
is negative we get the same equation \eqref{eq-omega},
thus we keep both solutions of \eqref{eq-omega}.

\noindent Replacing in \eqref{eq-omega-ex} the values
of $x$ and $y$ given in \eqref{eq:biharmonichelices2} we find
$$
\beta' -2\,m(c\,\sin \alpha_0 \sin \beta - d\,\sin \alpha_0 \cos \beta)
=\ell\,\cos\alpha_0 + 2\, m\, b\, \sin^2\alpha_0+\omega_{1,2}\,.
$$
Finally, taking into account the integral of the latter equation,
 the value of $z$ given in \eqref{eq:biharmonichelices2}
 becomes the desired ones.

The case of Type II and of Type III can be derived
in a similar way.
\end{proof}
\begin{remark}
We point out that the ODE \eqref{eq-diff-betaII} and \eqref{eq-diff-betaIII}
can be written as the Riccati equation with constant coefficients
\begin{equation}\label{eq-riccati}
x'=ax^2+bx+c,\quad a,b,c\in\r,
\end{equation}
which is a separable equation.

The ODE \eqref{eq-diff-beta} is of type
$$
x'=a \cos x + b \sin x +c,\quad a,b,c\in\r,
$$
which can be reduced to \eqref{eq-riccati}.
\end{remark}
\begin{remark}
The proper biharmonic curves of Type I, given in
\eqref{eq:biharmonichelices1}, lie in the ``round cylinder''
$$
S=\left\{ (x,y,z)\in M\,:\, (x-c)^2+(y-d)^2=b^2\,\sin^2\alpha_0\right\}
$$
and they are geodesics of $S$. The surface $S$ is invariant
by traslations along the $z$ axis, which are isometries
with respect to the Cartan-Vranceanu metric. Similarly, the
proper biharmonic curves of Type II and of Type III are
geodesics of the ``cylinders'' $y= x \tan\beta_0 +a$ and
$x=x_0$ respectively. For the case $\ell^2=4m,\,\ell\neq
0$, i.e. the case of the $3$-dimensional sphere, it was
proved, in \cite{RCSMCO1}, that the proper biharmonic
curves are geodesics on the Clifford Torus which is a
$SO(2)$-invariant surface of $\s^3$. We can conclude that
any proper biharmonic curve of the Cartan-Vranceanu spaces
is a geodesic on a surface which is invariant under the
action of an $1$-parameter  group of isometries.
\end{remark}
\noindent{\bf Acknowledgements}. The second author wishes
to thank the organizers of ``The Seventh International
Workshop on Differential Geometry
and its Applications - Deva - September 2005'' for their
exquisite hospitality and the opportunity of presenting a
lecture. The third author
wishes to thank the Dipartimento di Matematica e
Informatica,
Universit\`a di Cagliari, for hospitality during the
preparation of this paper.

\end{document}